%% file: r3k.tex
\documentclass[12pt]{article}

\usepackage{bm}

\makeatletter
\newfont{\footsc}{cmcsc10 at 8truept}
\newfont{\footbf}{cmbx10 at 8truept}
\newfont{\footrm}{cmr10 at 10truept}
  \makeatother
\title {\bf 
A step forwards on the Erd\H {o}s-S\'{o}s problem\\
concerning the Ramsey numbers $\bm{R(3,k)}$
}

\author{
Rujie Zhu\\[-.2ex]
\small College of Electrical Engineering\\[-0.8ex]
\small Guangxi University\\[-0.8ex]
\small Nanning, Guangxi 530004, China\\[-0.8ex]
\small \texttt{962637759@qq.com}\\[-1.2ex]
\\
Xiaodong Xu\\[-.2ex]
\small Guangxi Academy of Sciences\\[-0.8ex]
\small Nanning, Guangxi 530007, China\\[-0.8ex]
\small \texttt{xxdmaths@sina.com}\\[-1.4ex]
\\
and\\[-1.4ex]
\\
Stanis\l aw Radziszowski\\[-.2ex]
\small Department of Computer Science \\[-0.8ex]
\small Rochester Institute of Technology \\[-0.8ex]
\small Rochester, NY  14623, USA \\[-0.8ex]
\small \texttt{spr@cs.rit.edu}\\\
}

\begin{document}

\def\R{{\cal R}}
\newtheorem{thm}{Theorem}
\newtheorem{con}{Construction}
\newtheorem{col}{Corollary}
\newtheorem{conj}{Conjecture}

\maketitle

\begin{abstract}
Let $\Delta_s=R(K_3,K_s)-R(K_3,K_{s-1})$, where $R(G,H)$
is the Ramsey number of graphs $G$ and $H$ defined as the smallest
$n$ such that any edge coloring of $K_n$ with two colors contains
$G$ in the first color or $H$ in the second color. In 1980,
Erd\H {o}s and S\'{o}s posed some questions about the growth of
$\Delta_s$. The best known concrete bounds on $\Delta_s$ are
$3 \le \Delta_s \le s$, and they have not improved since the stating
of the problem. In this paper we present some constructions,
which imply in particular that $R(K_3,K_s) \ge R(K_3,K_{s-1}-e) + 4$. 
This does not improve the lower bound of 3 on $\Delta_s$,
but we still consider it a step towards to understanding
its growth. We discuss some related questions
and state two conjectures involving $\Delta_s$, including
the following:
for some constant $d$ and all $s$ it holds that
$\Delta_s - \Delta_{s+1} \leq d$. We also prove that
if the latter is true, then
$\lim_{s \rightarrow \infty} \Delta_s/s=0$.
\end{abstract}

{\bf AMS classification subjects:} 05C55

\eject


\bigskip
\section{Notation and Overview}

\bigskip
In this paper all graphs are simple and undirected.
The vertex set of graph $G$ is denoted by $V(G)$,
$n(G)=|V(G)|$, the edge set by $E(G)$,
and the set of neighbors of a vertex $v$ in $G$
will be written as $N_G(v)$.
The independence number of $G$, denoted $\alpha(G)$,
is the order of the largest independent set in $G$.
The graph induced in $G$ by the set of vertices $S \subset V(G)$
will be denoted by $G[S]$.
For $v \in V(G)$ and $e \in E(G)$, $G-v$ will stand
for $G[V\setminus \{v\}]$, and $G-e$ for the graph
$G$ with edge $e$ removed.

\medskip
For graphs $G$ and $H$, the {\em Ramsey number} $R(G,H)$ is the smallest
positive integer $n$ such that every coloring of the edges of $K_n$
with two colors contains a monochromatic $G$ in the first color or
a monochromatic $H$ in the second color. If the edges in the
first color are interpreted as a graph $F$ and those in the second
color as its complement, then $R(G,H)$ can be defined
equivalently as the smallest $n$ such that every $G$-free graph
on $n$ vertices contains $H$ in the complement. If $G=K_s$
and $H=K_t$ then we will write $R(s,t)$ for $R(G,H)$.
Any $G$-free graph $F$ on $n$ vertices without $H$
in the complement will be called a $(G,H;n)$-graph.
An $(s,t;n)$-graph will mean the same as a $(K_s,K_t;n)$-graph.
A regularly updated survey by the third author \cite{ds1}
lists the values and the best known bounds
on various types of Ramsey numbers.

In the sequel we will be concerned almost exclusively
with the Ramsey numbers $R(3,G)$ and $(3,G;m)$-graphs for
$G$ being $K_s$ or $K_s-e$. Observe that $R(3,G)=m+1$
if and only if $m$ is the largest integer such that there exists
a $(3,G;m)$-graph. Note also that in triangle-free
graphs the neighborhoods are independent sets.

\bigskip
The asymptotics of $R(3,s)$ was extensively
studied and now it is quite well understood. It is known that

$$\Big( 1/4+o(1)\Big) {s^2 \over {\log s}}
\leq R(3,s) \leq
\Big( 1+o(1)\Big) {s^2 \over {\log s}}.$$

\bigskip
\noindent
In 1995, Kim \cite{Kim} using probabilistic method improved
lower bound asymptotics to $R(3,s)=\Omega(s^2/\log s)$.
More detailed work followed, and finally
the lower bound constant 1/4 was obtained recently by Bohman
and Keevash \cite{BohK2}, and independently by Fiz Pontiveros,
Griffiths and Morris \cite{FizGM}. The upper bound constant 1 is
implicit in a 1983 paper by Shearer \cite{She1}, and it also
can be stated without $o(1)$ for $s\ge 3$ as

$$R(3,s+1) \leq
{{(s-1)^2} \over {\log s-1+s^{-1}}} +1.\eqno{(1)}$$

\bigskip
However,
the difference $R(3,G)-R(3,H)$ for concrete ``consecutive"
$G$ and $H$ is still very difficult to estimate, even starting
with rather small cases. In general, for
$K_s$ and $K_s-e$, all we know is the following:

\bigskip
\noindent
Easy old bounds \cite{BEFS}, see Section 2 and Construction 1 in Section 3,
$$3 \le R(3,K_s)-R(3,K_{s-1}) \le s,\hspace{1.3cm}\eqno{(2)}$$

\medskip
\noindent
trivial bounds implied by the monotonicity of Ramsey numbers
$$R(3,K_{s-1}) \le R(3,K_s-e) \le R(3,K_s),\eqno{(3)}$$

\medskip
\noindent
and a result obtained in this paper (Corollary 7 in Section 4)
$$4 \le R(3,K_{s+1})-R(3,K_s-e).\hspace{1.4cm}\eqno{(4)}$$

\bigskip
Many attempts were made to improve on some part of (2) or (3),
to no avail. We believe that our relatively simple constructions
proving inequality (4) in Section 4 form an interesting step
towards a better understanding of both (2) and (3). We pose it
as a challenge to improve over any of the inequalities in (2), (3)
or (4), or their combination as (4) combines parts of (2) and (3).


\section{Erd\H {o}s-S\'{o}s Problem}

\bigskip
\noindent
{\bf Problem.} Erd\H {o}s-S\'{o}s 1980 \cite{erd0,ChGr}

{\it
\noindent
Let $\Delta_s = R(3,s)-R(3,s-1)$. Is it true that}
$$\Delta_s \stackrel{s \rightarrow \infty}{\longrightarrow} \infty \  ? \ \ \ \ \
\big( \Delta_s/s\big) \stackrel{s \rightarrow \infty}{\longrightarrow} 0 \  ?\eqno{(5)}$$

\medskip
Only easy bounds on $\Delta_s$ as in (2) are known.
The upper bound $\Delta_s \le s$ is obvious since the maximum degree
of $(3,s)$-graphs is at most $s-1$. The lower bound $3 \le \Delta_s$
looks misleadingly simple, but it is not trivial (see Construction 1
in the next section).
It was argued in \cite{GoRa2} that a better understanding
of $\Delta_s$ may come from the study of $R(3,K_s-e)$
relative to $R(3,K_s)=R(3,s)$, since

$$\Delta_s = \Big( R(3,K_s)-R(3,K_s-e) \Big) +
\Big( R(3,K_s-e)-R(3,K_{s-1}) \Big).$$

\medskip
Recent progress on what we know for small cases is significant
\cite{GoRa1,GoRa2},
however some very simple-looking questions remain open.
For example, we do not even know whether
$R(3,K_s-e) - R(3,K_{s-1})$ is positive for all large $s$.
However, in Section 4 we prove (4), and in Section 5 we show
that the second part of (5) holds under the assumption that
there exists a constant $d$ for which
$\Delta_s - \Delta_{s+1} \le d$ for all $s$.


\input body.tex


\bigskip
\section*{Acknowledgments}

The work of the second author was partially supported by
the National Natural Science Foundation (11361244) and
the Guangxi Natural Science Foundation (2011GXNSFA018142).


\medskip

\input ref.tex
\end{document}

%% file: body.tex
\section{Previous Constructions}

\medskip
Burr, Erd\H{o}s, Faudree and Schelp \cite{BEFS}
in 1989 gave a general lower bound construction yielding
$R(k,s+1) \ge R(k,s)+2k-3$ for $k, s \ge 2$. For
$k=3$ it is equivalent to the following construction,
which implies $\Delta_s \ge 3$.

\bigskip
\noindent
{\bf Construction 1.} \cite{BEFS}

\noindent{\it
For $s\ge 2$, given any $(3,s;n)$-graph,
we can extend it to a $(3,s+1;n+3)$-graph.
}

\medskip
\noindent
{\bf Proof.} 
Let $u \in V(G)$ be any vertex of a $(3,s;n)$-graph $G$.
Note that since $\alpha(G)<s$ we have $\deg_G(u) < s$.
A $(3,s+1;n+3)$-graph $G'$ extending $G$ is defined on
the set of 3 more vertices $V(G')=V(G) \cup \{v,x,y\}$
with the set of edges
$E(G')=E(G)\cup \{vw \ |\  uw \in E(G)\} \cup \{ ux, xy, yv\}$.
Consider any independent set $I$ in $G'$, and the
cardinality $t$ of its intersection with $\{u,v,x,y\}$.
If $t\le 1$ then $|I| \le \alpha(G)+1$, otherwise $t=2$
and we must have that at least one of the vertices $u$
and $v$ is in $I$. Thus $|I \setminus \{u,v,x,y\}| < s-1$,
and hence $G'$ is a $(3,s+1;n+3)$-graph.
$\diamond$

\bigskip
Construction 1 was generalized in \cite{XXR} as described
in the next Theorem 2, which in turn implies the lower
bounds in the following Theorem 3 \cite{XSR}.

\bigskip
\noindent
{\bf Theorem 2.} \cite{XXR}
{\it
For every $k \ge 3$ and $s, t \ge 2$,
given any $(k,s)$-graph $G$ and $(k,t)$-graph $H$,
if both $G$ and $H$ contain an induced subgraph
isomorphic to some $K_{k-1}$-free graph $M$, then}
$$R(k,s + t - 1) \ge  n(G) + n(H) + n(M) + 1.$$

\bigskip
\noindent
{\bf Theorem 3.} \cite{XSR}

\noindent{\it
If \ $2 \leq s \leq t$ and $k \geq 3$, then
\begin{displaymath}
R(k,s+t-1) \ge R(k,s) + R(k,t) + \left\{ {\begin{array}{*{20}l}
{k - 3,\;\;\;\textit{if}\ s = 2;}\\
{k - 2,\;\;\;\textit{if}\ s \ge 3\ or\ k \ge 5.}\\
\end{array}} \right.
\end{displaymath}
}


\section{New Constructions}

\medskip
We present two simple constructions, the second one generalizing
the first, which together apparently add some new understanding
of (2) and (3) and how they imply (4).

\bigskip
\noindent
{\bf Construction 4.}
{\it
For $s, t \ge 3$, given any $(3,s+1;m)$-graph $G$
and $(3,t+1;n)$-graph $H$, we construct from $G$ and $H$
a $(3,s+t;m+n-2)$-graph $F$.
}

\medskip
\noindent
{\bf Proof.} 
Let $G$ be any $(3,s+1;m)$-graph and $H$ any $(3,t+1;n)$-graph,
on disjoint sets of vertices, and consider arbitrary two vertices
$u \in V(G)$ and $v \in V(H)$.
We will construct a $(3,s+t;m+n-2)$-graph
$F$ on the vertex set $V(F)=V(G) \cup V(H) \setminus \{u,v\}$.
Denote $X_G=N_G(u)$, $Y_G=V(G-u)\setminus X_G$,
$X_H=N_H(v)$ and $Y_H=V(H-v)\setminus X_H$, so
$V(F)$ is partitioned into $X_G \cup Y_G \cup X_H \cup Y_H$.
The set of edges of $F$ is defined by
$E(F)=E(G-u) \cup E(H-v) \cup \{ xy\ |\ x \in X_G, y \in X_H\}$.
Clearly, graph $F$ is triangle-free and it
has the right number of vertices. We need to
show that $\alpha(F)<s+t$. Note that $F$ contains a complete
bipartite graph with partite sets $X_G$ and $X_H$, and thus
for any independent set $I$ in $F$ we have
$x_G=|I \cap X_G|= 0$ or $x_H=|I \cap X_H|=0$, and
also it holds that
$y_G=|I \cap Y_G|\le s-1$ and $y_H=|I \cap Y_H|\le t-1$.
Since $|I \cap V(G)|\le s$ and $|I \cap V(H)|\le t$,
for $x_G>0$ we have $|I| \le s+y_H$, and for $x_H>0$
we have $|I| \le y_G+t$. In both cases $|I|\le s+t-1$
as required.
$\diamond$

\medskip
We observe that Construction 4 (with adjusted $s$ and $t$)
gives an alternate proof of a part of Theorem 3 for $k=3$,
but also $R(3,s+t) \geq R(3,s+1)+R(3,t+1)-3$.
In the next construction we exploit in more detail
the structure of the base graphs $G$ and $H$.

\bigskip
\noindent
{\bf Construction 5.}
{\it
For $s, t \ge 3$, given any $(3,s+1;m)$-graph $G$
which has two nonadjacent vertices with at most $c_G$
common neighbors, and any $(3,t+1;n)$-graph $H$
which has two nonadjacent vertices with at most $c_H$
common neighbors, we construct from $G$ and $H$ a
$(3,s+t-1;m+n-c_G-c_H-4)$-graph $F$.
}

\medskip
\noindent
{\bf Proof.} 
Let $G$ and $H$ be as stated above, where $u_1, u_2$
have $c_G$ common neighbors in $G$ and $v_1, v_2$
have $c_H$ common neighbors in $H$, respectively.
We partition the set $V(G)\setminus \{u_1,u_2\}$
into $X_G^{12} \cup X_G^1 \cup X_G^2 \cup Y_G$,
where $X_G^{12}=N_G(u_1) \cap N_G(u_2)$,
$X_G^1=N_G(u_1) \setminus N_G(u_2)$,
$X_G^2=N_G(u_2) \setminus N_G(u_1)$, and
$Y_G=\{u \in V(G)\ |\ uu_1 \notin E(G)$ and $uu_2 \notin E(G)\}$.
Similarly, we set the partition
$V(H)\setminus \{v_1,v_2\}=
X_H^{12} \cup X_H^1 \cup X_H^2 \cup Y_H$ by considering
four possible adjacencies to vertices $v_1, v_2$ in $H$.
Obviously, $|X_G^{12}|=c_G$ and $|X_H^{12}|=c_H$.
We will construct graph $F$ on the set of vertices
$X_G^1 \cup X_G^2 \cup Y_G \cup X_H^1 \cup X_H^2 \cup Y_H$,
which has cardinality as needed. The set of edges of $F$ is
defined by
$$E(F)=E(G[X_G^1 \cup X_G^2 \cup Y_G]) \cup
E(H[X_H^1 \cup X_H^2 \cup Y_H])
\cup K(X_G^1,X_H^1) \cup K(X_G^2,X_H^2),$$

\noindent
where $K(X_G^1,X_H^1)$ and $K(X_G^2,X_H^2)$ are the edges of two complete
bipartite graphs between indicated pairs of sets. It remains to be shown
that $\alpha(F) \le s+t-2$. Let $I \subset V(F)$ be any independent
set in $F$, and denote by $x_G^1, x_G^2, y_G, x_H^1, x_H^2, y_H$
the orders of intersection of $I$ with the corresponding parts
of $V(F)$. Similarly as in Construction 4, we have
($x_G^1=0$ or $x_H^1=0$) and ($x_G^2=0$ or $x_H^2=0$).
Furthermore, $x_G^1+y_G, x_G^2+y_G \le s-1$, $y_G \le s-2$, and
$x_H^1+y_H, x_H^2+y_H \le t-1$, $y_H \le t-2$.
If $x_G^1>0$ and $x_G^2>0$, then $x_H^1=x_H^2=0$
and thus $|I| \le \alpha(G)+y_H \le s+(t-2)$.
If $x_H^1>0$ and $x_H^2>0$, then $x_G^1=x_G^2=0$
and thus $|I| \le \alpha(H)+y_G \le t+(s-2)$.
If only one of $x_G^1, x_G^2$ is
positive, say $x_G^1>0$, then $x_G^2=x_H^1=0$ and
$|I| \le (x_G^1+y_G)+(x_H^2+y_H) \le (s-1)+(t-1)$.
$\diamond$

\bigskip
One can look at Construction 5 as lowering the independence
number of a union of $G$ and $H$ by 2, but at the cost of dropping
$c_G+c_H+4$ vertices. In the next three corollaries we show how in
some cases we can further assume that $c_G=c_H=0$. We will say that a
$(3,s)$-graph $G$ is {\it edge minimal} if deletion of any of its
edges increases $\alpha(G)$ to $s$, and it is {\it edge maximal}
if addition of any edge creates a triangle. A $(3,s)$-graph $G$ is
called {\it bicritical} if it is both edge minimal and maximal.

\bigskip
\noindent
{\bf Corollary 6.}
For $s \ge 3$ and $m=R(3,s+1)-1$,
if there exists a $(3,s+1;m)$-graph which is not bicritical,
then $\Delta_{s+2} \ge 4$.

\medskip
\noindent
{\bf Proof.}
Let $G'$ be a $(3,s+1;m)$-graph which is not bicritical. If it is
not edge minimal, then the removal of some edge $e=uv \in E(G)$
gives a $(3,s+1;m)$-graph $G=G'-e$, in which vertices $u$ and $v$
have no common neighbors. If $G'$ is not edge maximal, let $G=G'$.
In either case we have a $(3,s+1;m)$ graph $G$ with $c_G=0$
which will be used with Construction 5.

Let $H$ be the well known unique $(3,4;8)$-graph with 10 edges,
which is the cycle on 8 vertices $v_1v_2\cdots v_8$ with two
consecutive main diagonal edges $v_1v_5$ and $v_2v_6$.
The vertices $v_3$ and $v_7$ have no common neighbors.
We will use this $H$ with $n=8,\  t=3$ and $c_H=0$.
The graph $F$ resulting from $G$ and $H$ by applying
Construction 5 is a $(3,s+2;m+4)$-graph, which proves
the claim that $\Delta_{s+2} \ge 4$.
$\diamond$

\bigskip
Table 1 presents known values and bounds on $R(3,K_s)$, $R(3,K_s-e)$
collected in \cite{GoRa1,GoRa2} and $\Delta_s$ for $s\le 16$.
We note that for $s \le 9$, i.e. for which the exact value of
$R(3,s)$ is known, there exist non-bicritical $(3,s;R(3,s)-1)$-graphs
for $s \in \{4,6,7,8\}$.

\bigskip
\noindent
{\bf Corollary 7.}
$R(3,s+1) \ge R(3,K_s-e) + 4$,
for $s \ge 2$.

\medskip
\noindent
{\bf Proof.}
Let $m=R(3,K_s-e)-1$.
Observe that every $(3,K_s-e;m)$-graph is a $(3,K_s;m)$-graph
after removal of any of its edges, furthermore the endpoints
of the removed edge share no common neighbors, since otherwise
the original graph would have a triangle. Using Construction 5
for any such $(3,K_s)$-graph as $G$, and the $(3,4;8)$-graph $H$
as in the proof of Corollary 6, gives a
$(3,s+1;m+4)$-graph $F$, which proves the lower bound.
$\diamond$

\begin{table}
\begin{center}
\begin{tabular}{|c|c|c|c||c|c|c|c|}
\hline
$s$&$R(3,J_s)$&$R(3,K_s)$&$\ \ \Delta_s\ \ $&$s$&$R(3,J_s)$&$R(3,K_s)$&$\ \ \ \Delta_s\ \ \ $\cr
\hline
3&\ \ 5&\ \ 6&3&10&37&40--42&4--6\cr
4&\ \ 7&\ \ 9&3&11&42--45&47--50&\ 5--10\cr
5&11&14&5&12&47--53&52--59&\ 3--12\cr
6&17&18&4&13&55--62&59--68&\ 3--13\cr
7&21&23&5&14&59--71&66--77&\ 3--14\cr
8&25&28&5&15&69--80&73--87&\ 3--15\cr
9&31&36&8&16&73--91&82--98&\ 3--16\cr
\hline
\end{tabular}

\caption{$R(3,J_s)$ and $R(3,K_s)$,
$J_s=K_s-e$, for $s \le 16$ \cite{GoRa2}.}
\end{center}
\end{table}

\bigskip
In Table 1, for cases when only the bounds are given for
$R(3,K_s-e)$ ($J_s=K_s-e$) or $R(3,K_s)$, we believe that
the exact values are much closer to lower bounds than upper
bounds. Actually, we expect that in most open cases
the exact values are equal to the listed lower bounds.
The exceptions, if any, likely include some of the lower bounds
for $R(3,K_s-e)$, $s \in \{12,14,16\}$, which currently
are the only cases in the scope of Table 1 when
they are the same as the best known lower bounds
for $R(3,K_{s-1})$.

\bigskip
We end this section with one more corollary, which is a little
more general than Corollaries 6 and 7. Theorem 3 for $s \ge k=3$
gives the inequality $R(3,s+t-1) \geq R(3,s)+R(3,t)+1$.
The following Corollary 8 increases two terms of its right hand
side and decreases the constant only by 6.

\bigskip
\noindent
{\bf Corollary 8.}
$R(3,s+t-1) \geq R(3,K_{s+1}-e)+R(3,K_{t+1}-e)-5$, for $s,t \ge 3$.

\medskip
\noindent
{\bf Proof.}
Let $m=R(3,K_{s+1}-e)-1$ and $n=R(3,K_{t+1}-e)-1$.
Consider any $(3,K_{s+1}-e;m)$-graph $G'$ and any
$(3,K_{t+1}-e;n)$-graph $H'$. As in the proof
of Corollary 7, let $G=G'-e$ for some edge
$e \in E(G')$, then $G$ is a $(3,K_{s+1};m)$-graph which
has two nonadjacent vertices without common neighbors. Similarly,
obtain $(3,K_{t+1};m)$-graph $H$ from $H'$. Now, by applying
Construction 5 to graphs $G$ and $H$ we obtain graph $F$
witnessing the claimed lower bound.
$\diamond$

\bigskip

\section{Two Conjectures}

\bigskip
\noindent
Observe that
$$R(3,s+k) - R(3,s-1) = \sum_{i=0} ^{k} \Delta_{s+i}.\eqno{(6)}$$

\bigskip
We expect $\Delta_s$ to grow similarly as $s/\log s$ to account
for the asymptotics of $R(3,s)$ known to be $\Theta(s^2/\log s)$,
though with some small perturbations. $\Delta_s$
is actually known to be nonmonotonic as can be seen in Table 1
for $s$ between 4 and 6. However, we believe that such oscillations
are contained as stated in the following Conjecture 9, where we
anticipate that the decrease between consecutive $\Delta_s$ is bounded
by a constant.

Gy\'{a}rf\'{a}s, Seb\H{o} and Trotignon \cite{GST}
in their study of chromatic gaps, using
Theorems 2 and 3, showed that we can obtain lower bounds on
$R(3,s+k) - R(3,s)$ better than the obvious $3k$,
for $k\ge 2, s \ge 3$. In particular, we have $\Delta_s \ge 3$,
$\Delta_s+\Delta_{s+1}\ge 7$ and
$\Delta_s+\Delta_{s+1}+\Delta_{s+2}\ge 11$.

\bigskip
\noindent
{\bf Conjecture 9.}
{\it
There exists $d\ge 2$ such that for all $s \geq 2$
we have $\Delta_s - \Delta_{s+1} \leq d$.
}

\bigskip
Clearly, if $\Delta_s$ is nondecreasing for large $s$ then
$\lim _{s \rightarrow \infty} \Delta _s = \infty$, but
even if we could prove Conjecture 9 with $d=1$ for $s$ sufficiently
large (note that $\Delta_9 - \Delta_{10}\ge 2$),
it is not clear that it would help to prove
$\lim _{s \rightarrow \infty} \Delta _s = \infty$.
However, we will show that if Conjecture 9 is true then
it implies a positive solution to the second part
of the Erd\H{o}s-S\'{o}s problem.

\bigskip
\noindent
{\bf Theorem 10.}
{\it
If Conjecture $9$ is true, then
$\lim_{s \rightarrow \infty} \Delta_s/s=0$.
}

\bigskip
\noindent
{\bf Proof.}
For contradiction, suppose that Conjecture 9 holds, but there exists
$\epsilon>0$ such that $\Delta_s \geq \epsilon s$ for infinitely many $s$,
furthermore satisfying $s \ge 2d/\epsilon$.
Note that the latter implies $\epsilon s/d- 2\ge 0$.
Define $k=\lfloor \epsilon s/d \rfloor$, then observe
that $k \ge 2$, $k+1 > \epsilon s/d$
and $\epsilon s-kd\geq 0$. Now, assuming Conjecture 9, we have
$\Delta_{s+i}\ge \epsilon s - id$ for $0 \le i \le k$, and using (6) we
obtain the bound
$$R(3,s+k) - R(3,s-1) \geq (k+1)(\epsilon s -kd/2),$$
which gives
$$R(3,s+k)  > \epsilon^2s^2/2d.\eqno{(7)}$$

\bigskip
\noindent
On the other hand, the bound (1) with the above
constraints on $s,\epsilon$ and $d$ implies
\begin{eqnarray}
R(3,s+k) &\le& {{(s+k-2)^2} \over {\log (s+k-1)-1+(s+k-1)^{-1}}}+1 \nonumber \\
&<& {{(s+\epsilon s/d-2)^2} \over {\log (s+\epsilon s/d-2)-1}}+1 \nonumber \\
&<& s^2\left( {{(1+\epsilon /d)^2} \over {\log s-1}}+{1\over{s^2}}\right)
=s^2f(s,\epsilon,d), \nonumber
\end{eqnarray}

\medskip
\noindent
where for fixed $\epsilon$ and $d$ we have $\lim_{s\rightarrow\infty}f(s,\epsilon,d)=0$.
This contradicts inequality (7) for $s$ large enough, and
hence $\lim_{s\rightarrow \infty}\Delta_s/s=0$.
$\diamond$

\bigskip
While we expect that $\lim_{s \rightarrow \infty} \Delta_s =\infty$ is true, 
it can be very difficult to prove. Instead, we propose a weaker
statement in Conjecture 11, and we think that it might be provable by
constructive methods. This may be feasible by exploiting the techniques
used in asymptotic nonprobabilistic lower bound constructions for $R(3,s)$
such as those in \cite{chcd,CPR,KosPR}. So far such techniques are
weaker than the probabilistic methods, but they are more general
than the attempts of this paper.

\bigskip
\noindent
{\bf Conjecture 11.}
{\it There exists integer $k$ such that}
$$\lim _{s \rightarrow \infty} \sum_{i=0}^{k} \Delta_{s+i} =\infty.$$

\bigskip
Finally, we remark that the constructive methods studied in this paper
have applications beyond gaining new insights on the growth of $\Delta_s$,
like in the study of connectivity and hamiltonicity of Ramsey-critical
$(k,s;R(k,s)-1)$-graphs \cite{BP,GST,XSR}, chromatic gaps \cite{GST},
or in multicolor case, for Shannon capacity of graphs
with bounded independence number \cite{XR}.

\bigskip